\documentclass[a4paper]{article}  
\usepackage{amsmath, amstext, amsxtra, amsthm, amscd, amsgen, amsbsy, amsopn, amsfonts, latexsym, amssymb, amscd}

\begin{document}  
  


\def\proofend{\hbox to 1em{\hss}\hfill $\blacksquare
$\bigskip } 
\def\powser#1{\lbrack \lbrack #1 \rbrack \rbrack }  
  
\newtheorem{theorem}{Theorem}[section]  
\newtheorem{proposition}[theorem]{Proposition}  
\newtheorem{lemma}[theorem]{Lemma}  
\newtheorem{remark}[theorem]{Remark}  
\newtheorem{remarks}[theorem]{Remarks}  
\newtheorem{definition}[theorem]{Definition}  
\newtheorem{corollary}[theorem]{Corollary}  
\newtheorem{example}[theorem]{Example}  
\newtheorem{assumption}[theorem]{Assumption}  
\newtheorem{problem}[theorem]{Problem}  
\newtheorem{question}[theorem]{Question}  
\newtheorem{conjecture}[theorem]{Conjecture}
\newtheorem{rigiditytheorem}[theorem]{Rigidity Theorem}
\newtheorem{observation}[theorem]{Observation}

\def\Z{{\mathbb Z}}  
\def\R{{\mathbb R}}  
\def\Q{{\mathbb Q}}  
\def\C{{\mathbb C}}  
\def\N{{\mathbb N}}  
\def\H{{\mathbb H}}  
\def\Zp #1{{\mathbb Z }/#1{\mathbb Z}}  
\def\cpt{compact}  
\def\wt{wt}
\def\cowt{co\wt}
\def\achtel{\frac {\dim M} 8}

\def\torus{$T$}  
\def\mathtorus{T}  
\def\codim{{\rm{codim}\ }}
  
\def\b{bun\-dle}  
\def\pb{principal \b }  
\def\vb{vector \b }  
  
\def\mfd{manifold}
\def\LFF{Lefschetz fixed point formula}  
\def\isorank{symrank}  

\def\ell{\varphi _{ell}}
\def\eell{\widetilde \ell }  
\def\ddelta {{\widetilde \delta }}  
\def\eepsilon{{\widetilde \epsilon }}
\def\CC{C_0}
\def\oorder{o}   
\def\oha{\cal H}  
  
\def\paperref#1#2#3#4#5#6{\text{#1:} #2, {\em #3} {\bf#4} (#5)#6}  
\def\bookref#1#2#3#4#5#6{\text{#1:} {\em #2}, #3 #4 #5#6}  
\def\preprintref#1#2#3#4{\text{#1:} #2 #3 (#4)}

\hyphenation{man-i-fold equiv-a-riant in-te-ger mod-u-lo tor-sion
re-pre-sen-ta-tion di-men-sion-nal}


\title{Cyclic actions and elliptic genera}  
\author{Anand Dessai}
\date{}
\maketitle  
  
\begin{abstract}
\noindent
Let $M$ be a $Spin$-manifold with $S^1$-action and let $\sigma \in S^1$ be of finite order. We show 
that the indices of certain twisted Dirac operators vanish if the action of $\sigma $ has 
sufficiently large fixed point codimension. These indices occur in the Fourier expansion of the 
elliptic genus of $M$ in one of its cusps. As a by-product we obtain a new proof of a theorem of 
Hirzebruch and Slodowy on involutions. 
\end{abstract}  
  
\section{Introduction}\label{section intro}
Let $M$ be a smooth closed connected $Spin$-manifold with smooth $S^1$-action and let $\sigma 
\in S^1$ be the element of order two. Hirzebruch and Slodowy \cite{HiSl} showed that the elliptic 
genus of $M$ can be computed in terms of the transversal self-intersection of the fixed point 
manifold $M^\sigma $ and used this property to deduce a vanishing theorem for certain characteristic numbers which occur in the Fourier expansion of the elliptic genus of $M$ in 
one of its cusps. 

In this note we extend this vanishing theorem from involutions to cyclic 
actions of arbitrary order. Our main result (see Theorem \ref{theorem cyclic}) is used in \cite{De} 
to exhibit obstructions against the existence of positively curved metrics with symmetry on $Spin$-manifolds. The proof of Theorem \ref{theorem cyclic} relies on the rigidity theorem for the elliptic genus which we shall 
recall first. As a general reference for the theory of elliptic genera we recommend \cite{HiBeJu, 
La}. 

The elliptic genus $\Phi $, in the normalization considered in \cite{HiSl, Wi}, is a ring 
homomorphism from the oriented bordism ring to the ring of modular functions (with $\Zp 2$-character) 
for $\Gamma 
_0(2):=\{A\in SL_2(\Z )\; \mid \; A\equiv (\begin{smallmatrix} * & *\\0 & 
*\end{smallmatrix})\bmod 2\}$. In one of the cusps of $\Gamma _0(2)$ (the 
signature cusp) the Fourier expansion of $\Phi (M)$ has an interpretation as a series of twisted 
signatures 
$$sign (M,\bigotimes _{n=1}^\infty 
S_{q^n}TM \otimes 
\bigotimes 
_{n=1}^\infty \Lambda_{q^n}TM)=sign(M)+2\cdot sign(M,TM)\cdot q +\ldots .$$
 Here $sign(M,E)$ denotes the index of the signature 
operator twisted with the complexified \vb \ $E_\C $, $TM$ denotes the tangent bundle and $\Lambda 
_t=\sum \Lambda 
^i\cdot t^i$ (resp. $S_t=\sum S^i\cdot t^i$) denotes the exterior 
(resp. symmetric) power operation. 

Following Witten \cite{Wi} the series above is best thought of as the ``signature'' of the free loop 
space ${\cal L}M$ of $M$ formally localized at the manifold $M$ of constant loops. We denote the 
series of twisted signatures by $sign(q,{\cal L}M)$. 

The main feature of the elliptic genus is its rigidity under $S^1$-actions. This phenomenon was first 
explained by Witten \cite{Wi} using standard conjectures from quantum field theory and then shown 
rigorously by Taubes and Bott-Taubes in \cite{BoTa, Ta} (cf. also \cite{Hi2, Li}). 

If $S^1$ acts by isometries\footnote{This is the case after averaging a given Riemannian metric over 
the $S^1$-action.} on $M$ and if $E$ is a vector bundle associated to $TM$ then the signature 
operator twisted with the complexified \vb \ $E_\C $ refines to an $S^1$-equivariant operator. Its 
index is a virtual $S^1$-representation which we denote by $sign_{S^1}(M,E) 
\in R(S^1)$. In particular, the expansion of the elliptic genus in the signature cusp refines to a 
series of equivariant twisted signatures $sign_{S^1}(q,{\cal L}M)\in R(S^1)\lbrack \lbrack q\rbrack 
\rbrack $. 

\begin{theorem}[Rigidity theorem \cite{BoTa, Ta}]\label{rigidity theorem} Let $M$ be a closed manifold with $S^1$-action. If $M$ is $Spin$ then each equivariant twisted 
signature occurring as coefficient in the series $sign_{S^1}(q,{\cal L}M)$ is constant as a character 
of $S^1$.\proofend 
\end{theorem}

\noindent We use the rigidity theorem 
to study the action of cyclic subgroups of $S^1$. Our investigation is inspired by work of Hirzebruch 
and Slodowy \cite{HiSl} on elliptic genera and involutions. As a motivation we shall briefly recall relevant aspects of their work. 

Let $M$ be a $Spin$-manifold with $S^1$-action and let $\sigma 
\in S^1$ be of order two. By the rigidity theorem the expansion of the 
elliptic genus in the signature cusp is equal to the $S^1$-equivariant expansion evaluated at $\sigma 
\in S^1$, i.e. $sign(q,{\cal L}M)=sign_{S^1}(q,{\cal L}M)(\sigma )$. The latter can be computed via the 
\LFF \ \cite{AtSi} as a sum of local contributions $a_F$ at the connected components $F$ of the fixed point manifold 
$M^{\sigma }$. Hirzebruch and Slodowy showed that $a_F$ is equal to the expansion of the
elliptic genus (in the signature cusp) of the transversal self-intersection $F 
\circ F $ (cf. \cite{HiSl} for details):
\begin{equation}\label{formula elliptic}
sign (q,{\cal L}M)=sign _{S^1}(q,{\cal L}M)(\sigma ) \end{equation} $$=\sum 
_{F\subset M^\sigma }sign(q,{\cal L}(F 
\circ F ))=sign(q,{\cal L}(M^\sigma \circ 
M^\sigma ))$$
Note that, by taking constant terms, one obtains the 
classical formula $sign (M)=sign(M^\sigma 
\circ M^\sigma )$ for the ordinary signature which holds for the 
larger class of oriented manifolds (cf. \cite{Hi1,JaOs}). 

Formula (\ref{formula elliptic}) has two immediate consequences. If the codimension of $M^\sigma $, 
$\codim M^\sigma 
:=\min _{F\subset M^\sigma } \codim F$, is greater than half of the dimension 
of $M$ then the series $sign (q,{\cal L}M)$ vanishes identically. If the codimension of $M^\sigma $ 
is equal to half of the dimension of $M$ then all the twisted 
 signatures occurring as coefficients of $q^n$, $n>0$, in the series $sign (q,{\cal L}M)$ vanish,
  i.e. $sign (q,{\cal L}M)=sign(M)$. 

If the codimension of $M^\sigma $ is less than half of the dimension of $M$ then formula 
(\ref{formula elliptic}) still gives some information on the action of the involution $\sigma $. 
Namely it implies that certain twisted Dirac operators have vanishing index provided that the
codimension of $M^\sigma $ is sufficiently large. These indices are related to the elliptic genus in 
the following way. Recall that the $q$-series $sign (q,{\cal L}M)$ is the expansion of the elliptic 
genus $\Phi (M)$ in one of the cusps of $\Gamma 
_0(2)$. In a different cusp (the $\hat A$-cusp) the expansion 
of $\Phi (M)$ may be described (using a suitable change of cusps) by 
$$\Phi _0(M):=q^{-\dim M/8}\cdot \hat A(M,\bigotimes _{n=2m+1>0}\Lambda _{-q^n}TM \otimes 
\bigotimes
_{n=2m>0}S_{q^n}TM)$$  
$$=q^{-\dim M/8}\cdot (\hat A(M) -\hat A(M,TM)\cdot q +\hat A(M,\Lambda ^2{TM}+TM)\cdot
q^2+\ldots ).$$ Here $\hat A(M,E)$ is a characteristic number of the pair $(M,E)$ which, in the 
presence of a $Spin$-structure, is equal to the index of the Dirac operator twisted with the 
complexified 
\vb \ $E_\C $. We call the series above the expansion of $\Phi (M)$ in the $\hat 
A$-cusp. 

Note that $\Phi 
_0(M)$ and $sign(q,{\cal L}(M))$ are different expansions of the same 
modular function $\Phi (M)$ and determine each other. By formula (\ref 
{formula elliptic}) $\Phi _0(M)=\Phi _0(M^\sigma \circ M^\sigma )$ which implies the following generalization of the Atiyah-Hirzebruch vanishing 
theorem for the $\hat A$-genus \cite{AtHi}. 
\begin{theorem}[\cite{HiSl}]\label{theorem involution} Let $M$ be a $Spin$-manifold with 
$S^1$-action and let $\sigma \in S^1$ be of order two. If $\codim M^\sigma >4r$ then the expansion of 
the elliptic genus of $M$ in the $\hat A$-cusp has a pole of order less than $\achtel -r$.\proofend 
\end{theorem}

\noindent
The reasoning indicated above also leads to obstructions against the existence of $S^1$-actions on highly 
connected manifolds which might be of independent interest. 
\begin{theorem}\label{observation} Let $M$ be a $k$-connected $Spin$-manifold. Assume $k\geq 4r$. If $M$ admits a non-trivial $S^1$-action then the 
expansion of the 
 elliptic genus of $M$ in the 
$\hat A$-cusp has a pole of order less than $\achtel -r$. 
\end{theorem}

\bigskip
\noindent
Note that for $r>0$ the $Spin$-condition follows from the connectivity assumption. We remark that the 
conclusion of Theorem \ref{observation} also holds if $M$ is a connected $Spin$-manifold with 
non-trivial $S^1$-action and $H^{4*}(M;\Q )=0$ for $0<*\leq r $ (see Section \ref{highly connected} 
for a proof). 

The next result extends Theorem \ref{theorem involution} to finite cyclic actions of arbitrary order . 
 
\begin{theorem}\label{specialcase}
 Let $M$ be a $Spin$-manifold with 
$S^1$-action and let $\sigma \in S^1$ be of order $\oorder \geq 2 $. If $\codim M^\sigma 
> 2\oorder \cdot r$ then the expansion of the elliptic genus of $M$ in the 
$\hat A$-cusp has a pole of order less than $\achtel -r$. 
\end{theorem}

\noindent
The theorem follows from a more general result (see Theorem \ref{theorem cyclic} and the proof in Section \ref{section proof}). As indicated above the proof of Theorem \ref{theorem involution} given in \cite{HiSl} is specific to 
actions of order two. To deal with the general situation we consider the expansion of the equivariant 
elliptic genus in the $\hat A$-cusp and study the local contributions of the $S^1$-fixed point 
components using the rigidity theorem. We close this section with some consequences of Theorem 
\ref{specialcase}. 

\begin{corollary} Let $M$ be a $Spin$-manifold with 
$S^1$-action. 
\begin{enumerate} 
\item Let $\sigma \in S^1$ be of order $3$. If $\codim M^\sigma 
> 0$ then $\hat A(M)$ vanishes. If $\codim M^\sigma 
> 6$ then $\hat A(M)$ and $\hat A(M,TM)$ vanish. If $\sigma $ acts with isolated fixed points then $\Phi (M)$ vanishes identically.
\item Let $\sigma \in S^1$ be of order $4$. If $\codim M^\sigma 
> 0$ then $\hat A(M)$ vanishes. If $\codim M^\sigma 
> 8$ then $\hat A(M)$ and $\hat A(M,TM)$ vanish. If $\sigma $ acts with isolated fixed points then $\Phi (M)$ is equal to the signature of $M$.
\item Let $\sigma \in S^1$ be of order $\oorder <\frac {dim M}2$. If  $\sigma $ acts with isolated fixed points then $\hat A(M)$ and $\hat A(M,TM)$ vanish.\proofend
\end{enumerate}
\end{corollary}

\section{Cyclic actions}\label{section cyclic actions}
In this section we state the main result of this note. Let $M$ be a connected $S^1$-manifold and let 
$\oorder \geq 2$ be a natural number. At a connected component $Y$ of the fixed point manifold 
$M^{S^1}$ the tangent bundle $TM$ splits equivariantly as the direct sum of $TY$ and the normal 
bundle $\nu $. The latter splits (non-canonically) as a direct sum $\nu =\bigoplus_{k\neq 0} \nu _k$ 
corresponding to the irreducible real $2$-dimensional $S^1$-representations $e^{i\cdot 
\theta }\mapsto \left (\begin{smallmatrix} 
\cos k\theta & -\sin k\theta \\ 
\sin k\theta &\cos k\theta
\end{smallmatrix}\right )$, $k\neq 0$. We fix such a decomposition of $\nu $. 
For each $k\neq 0$ choose $\alpha 
_k\in \{\pm 1\}$ such that $\alpha _k k\equiv \tilde k \bmod \oorder $, 
$\tilde k\in 
\{0,\ldots , \lbrack 
\frac o 2\rbrack \}$.
 On each vector bundle 
$\nu _k$ introduce a complex structure such that $\lambda \in S^1$ acts on $\nu _k$ by scalar 
multiplication with $\lambda ^{\alpha _k k}$. The $\alpha _k k's$ (taken with multiplicities) are 
called the rotation numbers of the $S^1$-action at 
 $Y$. Finally define
$$m_\oorder (Y):=(\sum _k d_k
\cdot \tilde k )/\oorder \quad \text{ and }\quad m_\oorder :=\min _{Y} m_\oorder (Y),$$
where $d_k$ denotes the complex dimension of $\nu_k$ and $Y$ runs over the connected components of 
$M^{S^1}$ (to keep notation light we have suppressed the dependence of $\nu $, $\nu _k$, $d_k$ on 
$Y$). We are now in the position to state 
\begin{theorem}\label{theorem cyclic} Let $M$ be a $Spin$-manifold with 
$S^1$-action. If $m_\oorder > r$ then the expansion of the elliptic genus of $M$ in the $\hat A$-cusp 
has a pole of order less than $\achtel -r$. 
\end{theorem}

\bigskip
\noindent
If $\sigma \in S^1$ has order $\oorder 
=2$ then $\tilde k\in 
\{0,1 
\}$ and $4\cdot m_2 (Y)$ is the codimension of the connected component 
of $M^\sigma $ which contains $Y$. Thus $\codim M^\sigma 
\leq 4\cdot m_2$ and one recovers Theorem \ref{theorem 
involution}. In general if $\sigma \in S^1$ has order $\oorder 
\geq 2$ then $\codim M^\sigma \leq 2\oorder \cdot m_\oorder $ 
and one obtains Theorem \ref{specialcase}. Note that without the $Spin$ condition the conclusion of 
the theorem fails in general, e.g. for complex projective spaces of even complex dimension (see 
however Remark \ref{rigidity remark}). 
\section{Proof of Theorem \ref{theorem cyclic}}\label{section proof}
We may assume that the dimension of $M$ is divisible by $4$ and that the fixed point manifold 
$M^{S^1}$ is not empty since otherwise $M$ is rationally zero bordant by the \LFF \ \cite{AtSi} and $\Phi (M)$ vanishes. We may also assume that the $S^1$-action lifts to the 
$Spin$-structure (otherwise the action is odd which forces the elliptic genus to vanish, see for 
example \cite{HiSl}). We fix an $S^1$-equivariant Riemannian metric on $M$. The proof is divided into 
three steps. 

Step 1: We describe the equivariant elliptic genus at $M^{S^1}$. Consider the expansion of $\Phi (M)$ 
in the $\hat A$-cusp. Recall that the coefficients are indices of twisted Dirac operators associated 
to the $Spin$-structure. Since the $S^1$-action lifts to the $Spin$-structure each index refines to a 
virtual $S^1$-representation and the series refines to an element of $R(S^1)\lbrack q^{-\frac 1 
2}\rbrack \lbrack 
\lbrack q\rbrack 
\rbrack $ which we denote by $\Phi _{0,S^1}(M)$. Note that $sign_{S^1}(q,{\cal L}M)$ 
and $\Phi _{0,S^1}(M)$ are different expansions of the same function. Hence the rigidity of 
$sign_{S^1}(q,{\cal L}M)$ (see Theorem \ref{rigidity theorem}) is equivalent to the rigidity of $\Phi 
_{0,S^1}(M)$, i.e. each coefficient of the series 
 $\Phi _{0,S^1}(M)$ is constant as a character of $S^1$.

Let $\lambda _0\in S^1$ be a fixed topological generator. By the 
\LFF \ \cite{AtSi} the series $\Phi _{0,S^1}(M)(\lambda _0)\in \C \lbrack q^{-\frac 1 2}\rbrack \lbrack 
\lbrack q\rbrack \rbrack $
 is equal to a sum of local 
data 
$$\Phi_{0,S^1}(M)(\lambda _0)=\sum _Y \mu _Y(q,\lambda _0),$$
where $Y$ runs over the connected components of $M^{S^1}$. 

Recall from Section \ref{section cyclic actions} that we have decomposed the normal bundle $\nu $ of 
$Y$ as a direct sum $\bigoplus_{k\neq 0} \nu _k$ of complex vector bundles. Fix the orientation for 
$Y$ which is compatible with the orientation of $M$ and the complex structure of $\nu $. Let $\{\pm 
x_i\}$ denote the set of roots of $Y$ and let $\{x_{k,j}\}_{j=1,\ldots ,d_k}$ denote the set of roots 
of the complex vector bundle $\nu 
_k$. The local datum $\mu _Y(q,\lambda 
_0)$ may be described in cohomological terms as (cf. 
\cite{AtSi}, Section 3): 
\begin{equation}\label{local data}\mu _Y(q,\lambda _0)=\left \langle 
\prod _i \frac {x_i}
{f (q, x_i)}\cdot \prod _{k\neq 0\atop j=1,\ldots , d_k }\frac 1 {f (q, x_{k,j}+\alpha _k k\cdot 
z_0)}, \lbrack Y\rbrack \right \rangle \end{equation} Here $f (q,x)\in \C \lbrack 
\lbrack q^\frac 1 4\rbrack \rbrack \lbrack 
\lbrack x\rbrack \rbrack $ is equal to 
$$(e^{x/2}-e^{-x/2})\cdot q^{1/4}\cdot \frac {\prod _{n=2m>0}(1-q^n\cdot e^x)
\cdot (1-q^n\cdot e^{-x})}
{\prod _{n=2m+1>0}(1-q^n\cdot e^x)\cdot (1-q^n\cdot e^{-x})},$$ $\lambda _0=e^{z_0}$, $\lbrack 
Y\rbrack $ denotes the fundamental cycle of $Y$ and $\langle \quad ,\quad \rangle $ is the Kronecker 
pairing. In general each local datum $\mu 
_Y(q,\lambda _0)$ depends on $\lambda _0$. However, the sum
 $\sum _Y \mu _Y(q,\lambda _0)$ is equal to $\Phi_{0,S^1}(M)(\lambda _0)$ and therefore independent of 
$\lambda 
_0$ by the rigidity theorem. 

Step 2: Each local datum is the expansion of a meromorphic function on $\oha \times \C $ where $\oha 
$ denotes the upper half plane. As in the proof of the rigidity theorem given in \cite{BoTa} (cf. 
also \cite{DeJu, Hi2,Li}) modularity properties of these functions will be central for 
the argument. In this step we examine some of their properties. 

We begin to recall relevant properties of the series $f$ (see for example \cite{DeJu, HiBeJu}). For $0<\vert q\vert <1 $ and $z\in 
\C $ satisfying $\vert q\vert <\vert e^{z}\vert < \vert q\vert ^{-1}$ 
the series $f(q,z)$ converges normally to a holomorphic function. This function extends to a 
meromorphic function $\widetilde f(\tau ,z)$ on $\oha 
\times \C $ after the change of variables $q=e^{2\pi i\cdot 
\tau}$ where $\tau $ is in $\oha $. The function 
$\widetilde f(\tau ,z)$ is 
 elliptic in $z$ for the lattice $L:=4\pi i\cdot \Z 
 \langle 1,\tau \rangle $ and satisfies
$$\widetilde f(\tau 
,z+2\pi i)=- \widetilde f(\tau ,z), \widetilde f(\tau ,z+2\pi i\cdot 
\tau )=\widetilde f(\tau ,z)^{-1}, \widetilde f(\tau +2,z)=-\widetilde f(\tau ,z).$$ The zeros of $\widetilde f(\tau ,z)$ are simple and located at $L$ and $L+2\pi i$.

Let $q=e^{2\pi i\cdot \tau}$ and let $\lambda _0=e^{z_0}$ be a topological generator of $S^1$. In 
view of formula (\ref{local data}) and the properties of $f$ the local datum $\mu _Y(q,\lambda _0)$ 
converges to a meromorphic function $\widetilde \mu 
_Y$ on $\oha 
\times \C $ evaluated at $(\tau , z_0)$. We proceed to explain how this function is related to $\widetilde f$. For a function $F$ in the variables $x_i, x_{k,j}$ which is smooth in the origin let ${\cal 
T}(F)$ denote the Taylor expansion of $F$ with respect to $x_i, x_{k,j}=0$. It follows from formula 
(\ref{local data}) that $\widetilde \mu 
_Y$ is related to $\widetilde f$ by (see for example \cite{DeJu}):
$$\widetilde \mu _Y(\tau , z_0)=\left \langle {\cal T}
\left (\prod _i \frac {x_i} {\widetilde f(\tau , 
x_i)}\cdot \prod _{k\neq 0\atop j=1,\ldots ,d_k} \frac 1 {\widetilde f(\tau , x_{k,j}+\alpha _k 
k\cdot z_0)}\right ), \lbrack Y\rbrack \right \rangle $$ The properties of $\widetilde f$ stated 
above imply corresponding properties for $\widetilde \mu 
_Y$. In particular, $\widetilde \mu _Y$ is elliptic for 
the lattice $L$ and satisfies 
$$\widetilde 
\mu _Y(\tau +1,z)=(-1)^{\dim M/4}\cdot \widetilde 
\mu _Y(\tau ,z),\quad \widetilde 
\mu _Y(\tau ,z+2\pi i)=\pm \widetilde \mu _Y(\tau ,z).$$
For fixed $\tau \in \oha $ the poles of $\widetilde \mu _Y$ are contained in $\frac 1 n\cdot L$ for 
some $n\in 
\N $ depending on the rotation numbers of the $S^1$-action at $Y$ (see for example \cite{DeJu, HiBeJu}).

In general $\widetilde \mu _Y(\tau ,z)$ depends on $z$. If $\lambda 
=e^z$ is a topological generator of $S^1$, i.e. if $z/(2\pi i )$ is irrational, then $\Phi 
_{0,S^1}(M)(\lambda )$ converges to the sum $\sum 
_Y \widetilde \mu _Y(\tau ,z)$ by the \LFF \ and the latter is independent of $z$ by the rigidity theorem. Note that the original data may be recovered from $\widetilde \mu 
_Y(\tau ,z)$  by taking the 
expansion of $\widetilde \mu _Y(\tau ,z)$ with respect to $\tau 
\mapsto 
\tau +2$. 

Step 3: In the final step we study the series $\sum _Y \mu _Y$ in terms of the sum $\sum 
_Y \widetilde \mu _Y(\tau ,s(\tau ))$ where $s:\oha \to \C $ approximates $\tau \mapsto \frac 2 \oorder 
\cdot 2\pi i\cdot \tau $. We choose 
$s(\tau )$ in such a way that $\widetilde 
\mu _Y(\tau ,s(\tau ))$ is periodic with respect to $\tau \mapsto 
\tau +N$ for some $N\in \N $ (see below).

Note that in general the series $\mu 
_Y(q,\lambda )$ does not converge if $\lambda $ is close to $e^{\frac 
2 \oorder \cdot 2\pi i\cdot \tau }$ and the $q^{\frac 1 N}$-expansion of $\widetilde 
\mu _Y(\tau ,s(\tau ))$, denoted by $a_Y$, is different from the corresponding contribution $\mu _Y(q,\lambda _0)$ in the \LFF \ for $\Phi_{0,S^1}(M)(\lambda _0)$. 
In particular, we cannot compare $\mu 
_Y(q,e^{s(\tau )})$ and $\widetilde \mu _Y(\tau ,s(\tau ))$ directly. 
However, since the sum $\sum _Y \widetilde \mu _Y(\tau ,z)$ is independent of $z$ the sum $\sum _Y 
a_Y$ is equal to the elliptic genus in the $\hat A$-cusp (see last step). Using the properties of 
$\widetilde 
\mu _Y$ described above and the assumption on $m_\oorder $ we will show
 that $\sum _Y a_Y$ has a pole of order less than $\achtel -r$. This 
will complete the proof. 

Here are the details. The discussion in the last step implies that the poles of $\widetilde 
\mu 
_Y$, $Y\subset M^{S^1}$, are contained in $\frac 1 n \cdot L$ for some 
$n\in \N $. Choose $s(\tau ):=(1-\beta )\cdot \frac 2 {\oorder }\cdot 2\pi i \cdot \tau $, where 
$\beta $ is a fixed rational positive number $\ll \frac 1 n$. Hence, $s(\tau )$ is close to $\frac 2 
\oorder 
\cdot 2\pi i\cdot \tau $ and $\tau \mapsto \widetilde \mu _Y(\tau ,s(\tau ))$ is 
holomorphic on $\oha $ for every $Y$. Using $\alpha 
_k k\equiv \tilde k \bmod \oorder $, $\tilde k\in 
\{0,\ldots , \lbrack 
\frac o 2\rbrack \}$, and the transformation property 
$\widetilde f(\tau ,z+4\pi i\cdot\tau )=\widetilde f(\tau ,z)$ one computes that $\widetilde \mu 
_Y(\tau, s(\tau ))$ is (up to sign) equal to $\langle {\cal T}(A_Y), \lbrack Y\rbrack \rangle $, 
where 
$$A_Y:=\prod _i \frac {x_i} {\widetilde f(\tau , 
x_i)}\cdot \prod _{k\neq 0\atop j=1,\ldots ,d_k} \frac 1 {\widetilde f(\tau , x_{k,j}+2\cdot (\frac 
{\widetilde k}{\oorder}\cdot (1-\beta )-\beta _k) 
\cdot (2 \pi i \cdot \tau ) )}$$
and $\beta _k:=\beta \cdot \frac {\alpha 
_kk - \tilde k }\oorder $. 

Note that for some $N\in \N $ (depending on $\beta $ and the rotation numbers) every summand 
$\widetilde \mu _Y(\tau, s(\tau ))$ is periodic with respect to $\tau \mapsto \tau +N$. We claim that 
its expansion $a_Y\in \C 
\lbrack q^{-\frac 1 N}\rbrack \lbrack \lbrack q^{\frac 1 N}\rbrack 
\rbrack $ has a pole of order less than $\achtel 
-r$. 

Since the expansion of ${\cal T}\left (x_i/\widetilde f(\tau , x_i)\right )$ (with respect to $\tau 
\mapsto \tau +4$) is equal to $x_i/f(q, x_i)$ the expansion of 
\begin{equation}\tag{$\ast$}{\cal T}\left (\frac 1 {\widetilde f(\tau , x_{k,j}+2\cdot 
(\frac {\widetilde k}{\oorder}\cdot (1-\beta )-\beta _k)\cdot 
 (2 \pi i \cdot \tau ) )}\right )\end{equation}
can be easily computed in terms of $f$. The computation shows that the expansion of ($\ast$) has a 
pole of order $\leq \frac 1 4 - \frac {\widetilde k}{\oorder}\cdot (1-\beta )+\beta _k$. Since 
$m_\oorder (Y)\geq m_\oorder >r$ and $\beta $, $\beta _k$ are arbitrarily small it follows 
 that $a_Y\in \C \lbrack q^{-\frac 1 N}\rbrack \lbrack \lbrack 
q^{\frac 1 N}\rbrack \rbrack $ has a pole of order less than $\achtel- r$. As explained in the 
beginning of this step the sum $\sum _Y a_Y$ is equal to the expansion of the elliptic genus in the 
$\hat A$-cusp. Hence, $\Phi_{0}(M)\in 
\C 
\lbrack q^{-\frac 1 2}\rbrack\lbrack \lbrack q\rbrack \rbrack $ has a pole of order less than 
$\achtel- r$. This completes the proof.\proofend 

\begin{remark}\label{rigidity remark} Essentially the same reasoning applies to 
orientable $S^1$-manifolds (not necessarily $Spin$) for which the equivariant elliptic genus is rigid. The rigidity theorem is known 
 to hold for oriented manifolds 
 with finite second
  homotopy group \cite{Her} and for $Spin^c$-manifolds with first Chern class a torsion 
  class \cite{Despinc}. Theorem \ref{theorem cyclic} is also true for these manifolds.
\end{remark}

\section{Highly connected ${\mathbf {S^1}}$-manifolds}\label{highly connected}
In this section we adapt the arguments of \cite{HiSl} to study the elliptic genus of certain 
$S^1$-manifolds including highly connected manifolds. To begin with we recall the \LFF \ for twisted 
signatures. Let $M$ be an oriented closed $S^1$-manifold, $E$ an $S^1$-equivariant \vb \ over $M$ and 
$\sigma \in S^1$ the element of order $2$ . In the following we shall always assume that the fixed 
point manifold $M^\sigma $ is orientable (this is the case if $M$ is $Spin$ \cite{BoTa}). By the \LFF 
\ the equivariant twisted signature $sign _{S^1}(M,E)\in R(S^1)$ evaluated at $\sigma$ is equal to a sum of 
local data $a_{F,E}$ at the connected components $F$ of the fixed point manifold $M^\sigma $ 
$$sign_{S^1}(M,E)(\sigma )=\sum _F a_{F,E}.$$
The local contributions are given by (cf. \cite{HiSl})
$$a_{F,E}=\left \langle A_{F,E}, \lbrack F\rbrack \right \rangle $$
where
$$A_{F,E}=\prod _i\left (x_i\cdot \frac {1+e^{-x_i}}{1-e^{-x_i}}\right )\cdot \prod _j\left (y_j\cdot\frac {1+e^{-y_j}}{1-e^{-y_j}}\right )^{-1}
\cdot ch (E_{\vert F})(\sigma )\cdot e(\nu _F).$$
Here $\pm x_i$ (resp. $\pm y_j$) denote the formal roots of $F$ (resp. the normal bundle $\nu _F$ of 
$F$) for compatible orientations of $F$ and $\nu _F$, $e(\nu _F)$ is the Euler class of $\nu _F$ and 
$ch(E_{\vert F})$ denotes the equivariant Chern character of $E_{\vert F}$. The local datum $a_{F,E}$ is obtained by 
evaluating the cohomology class $A_{F,E}$ on the fundamental cycle $\lbrack F\rbrack $ via the 
Kronecker pairing $\langle 
\quad ,\quad \rangle $. Note that $a_{F,E}$ vanishes if $e(\nu _F)$ is a torsion class. Hence, the 
following lemma is immediate. 

\begin{lemma}\label{first lemma} Let $M$ and $E$ be as above and let $F\subset M^\sigma$ be 
of codimension $k$. If $H^{k}(F ;\Q )=0$ then the local datum $a_{F,E}$ vanishes.\proofend
\end{lemma}

\noindent
For the proof of the next lemma recall that the Euler class of the normal bundle of 
$i:F\hookrightarrow M$ is equal to $i^*(i_!(1))$, where $i_!:H^*(F;\Z )\to H^{*+k}(M;\Z )$ denotes 
the push forward in cohomology for the oriented normal bundle $\nu _F$. 

\begin{lemma}\label{second lemma} Let $M$ and $E$ be as above. If $H^{k}(M ;\Q )=0$ then $a_{F,E}$
 vanishes for any connected component $F\subset M^\sigma$ of codimension $k$.\proofend
\end{lemma}

\noindent
We shall now apply these observations to the elliptic genus. 
\begin{theorem} Let $M$ be a $Spin$-manifold. Assume that $H^{4*}(M;\Q )=0$ for $0<*\leq r $.
 If $M$ admits a non-trivial $S^1$-action then the 
expansion of $\Phi(M)$ in the $\hat A$-cusp has a pole of order less than $\achtel -r$. 
\end{theorem}

\bigskip
\noindent
{\bf Proof:} Let $\sigma \in S^1$ denote the element of order two. Arguing as in the proof of Theorem 
\ref{theorem cyclic} we may assume that the dimension of $M$ and the dimension of each connected 
component $F\subset M^\sigma $ is divisible by $4$. Consider the expansion $sign _{S^1}(q,{\cal L}M)$ 
of the $S^1$-equivariant elliptic genus in the signature cusp. By the rigidity theorem $sign 
_{S^1}(q,{\cal L}M)(\sigma )$ is equal to the non-equivariant expansion $sign (q,{\cal L}M)$. By the 
\LFF \ $sign _{S^1}(q,{\cal L}M)(\sigma )$ is a 
sum of local contributions $a_F$ at the connected components $F$ of $M^\sigma $:
$$sign (q,{\cal L}M)=sign  _{S^1}(q,{\cal L}M)(\sigma )=\sum _F a_F.$$
Note that each coefficient of the $q$-power series $a_F$ is the local contribution in the 
\LFF \ of an equivariant twisted signature evaluated at $\sigma \in S^1$. Since $H^{4*}(M;\Q )=0$ for 
$0<*\leq r $ the contribution $a_F$ vanishes if $\codim F\leq 4r$ (see Lemma \ref{second lemma}). If 
$\codim F> 4r$ then $a_F$ is equal to $sign (q,{\cal L}(F\circ F))$ (see formula (\ref{formula 
elliptic})). Hence, 
$$sign (q,{\cal L}M)=\sum _{\codim F> 4r} a_F=\sum _{\codim F\circ F>8r}sign  (q,{\cal L}(F\circ F)).$$
This implies that the expansion of $\Phi(M)$ in the $\hat A$-cusp has a pole of order less than 
$\achtel -r$.\proofend 

\bigskip
\noindent
Finally note that Theorem \ref {observation} is a direct consequence of the theorem above.

\bigskip
\noindent
Anand Dessai\\ e-mail: dessai@math.uni-augsburg.de\\ 
http://www.math.uni-augsburg.de/geo/dessai/homepage.html\\ Department of Mathematics, University of 
Augsburg, D-86159 Augsburg\\

\end{document}